\documentclass[11pt]{amsart}
\usepackage{amsmath,amssymb, amscd}
\usepackage{graphicx}

\usepackage{verbatim}

\def\sqr#1#2{{\vcenter{\hrule height.#2pt
        \hbox{\vrule width.#2pt height#1pt \kern#1pt
                \vrule width.#2pt}
        \hrule height.#2pt}}}

\newtheorem{theorem}{Theorem}[section]

\newtheorem{proposition}[theorem]{Proposition}

\theoremstyle{definition}
\newtheorem{definition}[theorem]{Definition}
\newtheorem{example}[theorem]{Example}

\newtheorem{question}[theorem]{Question}

\newtheorem{remark}[theorem]{Remark}

\newtheorem{discussion}[theorem]{Discussion}

\DeclareMathOperator{\n}{\mathbf n}
\DeclareMathOperator{\m}{\mathbf m}

\DeclareMathOperator{\Spec}{Spec}

\DeclareMathOperator\grad{grad}

\DeclareMathOperator\ord{ord}

\DeclareMathOperator\Proj{Proj}

\DeclareMathOperator{\tr.deg}{tr.deg}

\DeclareMathOperator{\N}{\mathbb N}

\DeclareMathOperator{\Rees}{Rees }
\DeclareMathOperator{\Char}{char  }

\DeclareMathOperator{\info}{info  }

\def\alert#1{\smallskip{\hskip\parindent\vrule%
\vbox{\advance\hsize-2\parindent\hrule\smallskip\parindent.4\parindent%
\narrower\noindent#1\smallskip\hrule}\vrule\hfill}\smallskip}

\begin{document}

\title[Abhyankar's Work on Dicritical Divisors ]
{Abhyankar's Work on Dicritical Divisors }


\author{William Heinzer}
\address{Department of Mathematics, Purdue University, West
Lafayette, Indiana 47907 U.S.A.}
\email{heinzer@purdue.edu}

\author{David Shannon}
\address{Department of Mathematics, Transylvania University, Lexington, KY 40508 U.S.A.}
\email{dshannon@transy.edu}

%\thanks{ }

\date \today


\maketitle
\bigskip

\begin{abstract}
We discuss the work of Abhyankar on dicritical divisors with a special focus on the
algebraic aspects of this work.   We also discuss related work on local quadratic transforms,
infinitely near points and Rees valuation rings of an ideal.

\end{abstract}

\baselineskip 19 pt

\section{Introduction} \label{c1}

Early in his career, in the context of working on the problem of the resolution of singularities,  Abhyankar
published in \cite[1956]{Ab1}  a paper that has turned out to be one of his most cited papers, ``On the
Valuations Centered in a Local Domain''.    In this paper he proves a theorem (Proposition~3 of \cite{Ab1})
that characterizes prime divisors of a  regular local domain.\footnote{For the
definition of prime divisor, see Definition~\ref{2.3}.    In Zariski-Samuel    \cite[page~95]{ZS2},   valuation domains with this
property are said to be prime divisors of the second kind.}
The characterization may be described as follows.

\begin{theorem} \label{1.20}
Let $R$ be an $n$-dimensional regular local domain with maximal ideal $M(R)  = M$ and  assume that $n \ge 2$.
Let $V$ be a prime divisor of
$R$ with  center  $M$ in $R$.\footnote{A valuation domain $V$ on the quotient field of $R$ that contains $R$
 is said to have {\it center} $M$ in $R$
 if the maximal ideal $M(V)$ of $V$ intersects with
 $R$ in $M$.  In this case we say $V$ {\it dominates} $R$.}   There exists a unique  finite sequence
\begin{equation} \label{1}
R ~ =  ~R_0  ~\subset ~ R_1 ~ \subset \cdots \subset ~ R_{h}  ~\subset ~ R_{h+1} ~= ~ V
\end{equation}
of regular local rings $R_j$, where  $\dim R_h \ge 2$ and    $R_{j+1}$ is the first local quadratic transform  of
$R_{j}$ along $V$   for each $j \in \{0, \ldots, h \}$,
and $\ord R_{h}   = V$.\footnote{For quadratic transforms, see for example \cite[pp. 569- 577]{Ab8} and \cite[p. 367]{ZS2}. The powers of the maximal ideal of a regular local domain $S$ define a rank one
discrete valuation domain  denoted  $\ord S$.   If $\dim S = d$,  then the residue field of $\ord S$ is a pure transcendental
extension of the residue field of $S$ of transcendence degree $d-1$.}
\end{theorem}

It follows from Theorem~\ref{1.20}  that  the residue field $V/M(V)$ of $V$ is a pure
transcendental extension of
the field $R_h/M(R_h)$ of transcendence degree   one less than $\dim R_{h} $.
 Therefore the residue field of $V$
is ruled as an extension field of the residue field of $R$.\footnote{A field extension
 $F \subset L$ is said to be {\it ruled}
if $L$ is a simple transcendental extension of  a subfield $K$ such that  $F \subset K$.}

The association of the prime divisor $V$ with the regular local ring $R_h$ in Equation~\ref{1}, and the
uniqueness of the sequence in Equation~\ref{1}
establishes a one-to-one correspondence between the prime divisors $V$ dominating the
regular local ring $R$ and the regular local rings $S$ of dimension at least 2
 that  dominate $R$ and are
obtained from $R$  by a finite sequence of local quadratic transforms    as in Equation~\ref{1}.
The regular local rings $R_j$ with $j \le h$ displayed in Equation~\ref{1} are
the {\bf infinitely near points } to $R$
along $V$.  In general, a regular local ring $S$ of dimension at least 2  is
called an {\bf infinitely near point } to $R$ if
there exists a sequence
$$
R ~ =  ~R_0  ~\subset ~ R_1 ~ \subset \cdots ~ \subset ~ R_{h}~ = ~ S,  \quad h ~\ge ~ 0
$$
of regular local rings $R_j$ of dimension at least 2,   where   $R_{j+1}$ is the first
local quadratic transform  of
$R_{j}$   for each $j$ with $0 \le j \le h-1$.
\cite[Definition~1.6]{L}.

 The Zariski-Abhyankar Factorization Theorem \cite[Theorem~3]{Ab1} implies that if $\dim R = 2$,
then every 2-dimensional regular local ring $S$ that birationally
dominates\footnote{A quasilocal extension  domain $S$    {\it birationally dominates}  $R$
if $S$ is contained in the quotient field of $R$ and $M(S) \cap R = M$. }
 $R$ is an infinitely near point
to $R$.     We record  in Theorem~\ref{1.1} implications of \cite[Theorem~3]{Ab1}.

\begin{theorem} \label{1.1}Let $R \hookrightarrow S$ be a birational extension of
2-dimensional  regular local domains.
\begin{enumerate}
\item If $R \ne S$, then $M(R)S$ is a proper principal ideal of $S$.
Therefore $S$ dominates a unique local quadratic transform $R_1$ of $R$.
\item There exists for some positive integer $\nu$ a sequence
$$
R~= ~ R_0 ~\subset ~R_1 ~\subset ~ \cdots ~ \subset R_{\nu} ~= ~ S,
$$
where $R_i$ is a local quadratic transform  of $R_{i-1}$ for each $i \in \{1, \ldots, \nu\}$.
The rings $R_i$ are precisely the regular local domains  that are subrings of $S$ and contain
$R$.

\end{enumerate}

\end{theorem}

 In the case where $R$ is a regular local domain with
  $\dim R \ge 3$,  there are many regular local rings $S$ that birationally dominate $R$
 with $\dim R = \dim S$  such that $S$ is not an infinitely near point to $R$.
 Perhaps the simplest examples are provided by local monoidal transforms, cf.   \cite{S}.

 \begin{example} \label{1.6}   Let $(R,\m)$ be a 3-dimensional regular local ring with maximal ideal
$\m = (x,y,z)R$, and let $S = R[\frac{y}{x}]_{(x, \frac{y}{x}, z)R[\frac{y}{x}]}$.
  Then $S$ is a 3-dimensional regular local ring
  that birationally dominates $R$ and $\m S = (x, z)S$ is a prime ideal of $S$ of height 2.  It follows that
   $S$ does not
  dominate a local quadratic transform of $R$.  Therefore  $S$  is not an infinitely near point of $R$.
Thus $S$ is a 3-dimensional regular local domain  that birationally dominates $R$,
but $S$ is not an infinitely near point to $R$.

Consider the blowup $\Proj R[\m t]$ of the maximal ideal $\m$ of $R$.  In the notation of Abhyankar,
$\Proj R[\m t]$  is the modelic blowup of $R$ at $\m$ and is denoted  $ \frak W(R, \m)$
in Definition~\ref{3.1}.
 Let $V$ denote the order valuation ring
of  $S$.  The center of $V$ on  the
modelic blowup  $\frak W(R, \m)$
  is the maximal ideal of the  2-dimensional regular local ring
$$
T ~ :=  ~R[x, \frac{y}{x}, \frac{z}{x}]_{(x, \frac{y}{x})R[x, \frac{y}{x}, \frac{z}{x}]}  ~ =   ~
R[ \frac{x}{z}, \frac{y}{z}, z]_{(z, \frac{y}{z})R[ \frac{x}{z}, \frac{y}{z}, z]}.
$$
The regular local ring $T$ is infinitely near to $R$;  indeed,  $T$ is a point in the first neighborhood of $R$.
However,  if $J$ is an     $\m$-primary ideal of $R$ that has $T$ as a base point,\footnote{See Section~5 for the 
definitions of base point  and finitely supported.}
then $J$ is not finitely supported \cite[Corollary~1.22]{L}.  Thus an  $\m$-primary ideal of $R$ such as
$J = (x^2, y, z^2)R$ is not finitely supported.
\end{example}

  Hence  in the case where $R$ is a regular local domain of dimension at least three,
  the one-to-one correspondence between prime divisors birationally
  dominating $R$ and regular local rings infinitely near to $R$  fails to include many of the regular local rings
  that birationally dominate $R$.  Indeed,  a prime divisor $V$ birationally dominating
 a $3$-dimensional regular local ring $R$ may be such that there exist infinitely  many $3$-dimensional regular
 local rings that birationally dominate $R$ and have $V$ as their order valuation ring, cf.
  \cite[Lemma~4.2 and Corollary~4.5]{Sa} and \cite[Example~2.6]{HK}.
  There can be, however,  no proper inclusion relations among  these $3$-dimensional regular local rings
  that have the same order valuation because of Theorem~\ref{sally1.4}, a result of  Sally
  \cite[Corollary~2.6]{Sally}.   Theorem~\ref{sally1.4}
   extends to higher dimensional regular local rings a result
  that is true for $2$-dimensional regular local rings
  by the Zariski-Abhyankar factorization theorem.

  \begin{theorem}  \label{sally1.4}    Let $R  \hookrightarrow T $ be a birational extension of
  $d$-dimensional regular local domain,  and let $V = \ord_R$ denote the rank one
  discrete valuation domain defined by the powers of the maximal ideal of $R$.  If $V$ dominates $T$,
  then $R = T$.
  \end{theorem}

During the early 1970's Abhyankar worked on the Jacobian problem that conjectures:  for polynomials $f_1, \ldots, f_n$
in the polynomial ring $k[x_1, \ldots, x_n]$ over a field $k$ with $\Char k = 0$,  if the determinant of the
Jacobian matrix $(\partial f_i/\partial x_j)$  is a nonzero constant, then
$k[f_1, \ldots, f_n] = k[x_1, \ldots, x_n]$. Abhyankar  returned to the Jacobian problem in the early 2000's.
    Using techniques involving characteristic sequences,  approximate roots and
Newton polygons (techniques described carefully in \cite{Ab4.5}),  Abhyankar developed solutions to the two variable
Jacobian problem
for various cases, the most general case being what he called the characteristic pair $2 + \epsilon$ case.  These results are
described in detail in the three papers  \cite{Ab8.1}, \cite{Ab8.2}, \cite{Ab8.3}.

Then in June and July  of 2008 in visits to Spain and France,  Abhyankar was introduced by Artal and Bodin
to some topological methods for analyzing the Jacobian problem, namely methods focused on the concept of
dicritical divisors.  Based on these initial conversations,  he was inspired to algebracize dicritical divisors.  He hoped
that this would lead to a fresh and powerful method to attack  the Jacobian problem,  as well as other related problems,
even problems  in mixed characteristic.

In a series of papers,  some by himself and some in collaboration with others,  Abhyankar created a precise and
general algebraic theory of dicritical divisors.  A key component in this development is the characterization of  prime divisors
obtained in his paper \cite[1956]{Ab1}.   Other fundamental components  include work Abhyankar had
done on: resolution of singularities, quadratic transformations, characteristic sequences and Newton polygons.
The classical work of Zariski on the theory of complete ideals in 2-dimensional regular local domains and related work
in ideal theory of Northcott and Rees were   also used.

\section{What is a dicritical divisor}

  In Definition~\ref{2.3}, we define   prime divisors of a
regular local ring and
the dicritical divisors of
the nonzero elements in the quotient field of a two-dimensional regular local domain.

\begin{definition} \label{2.3}
Let $R$ be a regular local ring with maximal ideal $M(R) = M$ and quotient field $L$.
\begin{enumerate}
\item
 A valuation domain $V$  with
quotient field $L$ that dominates $R$ is said to be a {\bf prime divisor} of $R$ if the residue field $V/M(V)$ has transcendence
degree $\dim R -1$ over the field $R/M$.
  Let  $D(R)^{\Delta}$   denote the set of prime divisors of $R$.

\item   Assume $\dim R = 2$, and let $z$ be a  nonzero element in $L$.
A prime divisor $V \in D(R)^{\Delta}$ is a {\bf dicritical divisor } of $z$ in $R$ if the image of $z$ in $V/M(V)$ is
transcendental over the field $R/M$.
Let $\mathfrak D(R,z)$ denote the set of dicritical divisors of $z$ in $R$.

\item
With $z$ and $R$ as in item 2,  we say that $z$ {\bf generates a special pencil} if there exists $x \in M \setminus M^2$ such
that $x^mz \in R$ for some $m \in \N$.
\end{enumerate}
\end{definition}

\begin{remark}
Let $R$ be a  2-dimensional regular local domain  with maximal ideal $M$ and quotient field $L$,
and let $ z$ be a nonzero element of $L$.
If $z$ or $1/z$ is in $R$, then the set $\mathfrak D(R,z)$  of dicritical divisors of $z$ in $R$ is
empty. Assume this does not hold and write $z = a/b$, where $a$ and $b$ are in
$R$ and have no common factors. Then $J = (a,b)R$ is $M(R)$-primary. In this
situation the dicritical divisors of $z$ are precisely the Rees valuations
of the ideal $J$ as defined for example in  Swanson and Huneke \cite[pages 187-210]{SH}.
\end{remark}

The set $\mathfrak  D(R,z)$ of dicritical divisors of $z$ in $R$ is a finite
set. An easy way to describe this set is to consider the extension  ring $R[a/b]$.
Let $t$ be an indeterminate over $R$ and consider the surjective $R$-homomorphism $R[t] \to R[a/b]$ defined
by $t \mapsto a/b$.  The kernel of this homomorphism is the prime ideal $(bt -a)R[t]$ and is contained in the
extension  $MR[t]$ of the maximal ideal $M$ of $R$ to $R[t]$.    It follows that  $MR[a/b]$ is a prime ideal of
height one,  and by the Krull-Akizuki Theorem,  the integral closure of $R[a/b]_{MR[a/b]}$ is a semilocal
PID. The dicritical divisors of $J = (a,b)R$ are precisely the DVRs
that birationally dominate $R[a/b]_{MR[a/b]}$.

 \begin{example}  \label{1.3}  Let $M(R) = M = (x,y)R$ and let $z = \frac{y^2}{x^3}$ and
$J = (x^3,y^2)R$. Consider the local ring
$$
S~~ = ~~R[y^2/x^3]_{MR[y^2/x^3]}~~ = ~~ \frac{R(t)}{(x^3t - y^2)R(t)}.
$$
Here $t$ is an indeterminate and $R(t)$ denotes the polynomial
ring $R[t]$ localized at the prime ideal $MR[t]$.  The ring $S$ is a
one-dimensional local domain. It is not integrally closed.  The element $\frac{y}{x}$ is
integral over $S$.
 This is related to the
fact that the ideal $I = (x^3, y^2, x^2y)R$ is integral over the ideal $J = (x^3, y^2)R$, or,  in
the terminology of Northcott-Rees, the ideal $J$ is  a reduction\footnote{An ideal $J$ is a {\it reduction} of
an ideal $I$ if $J$ is contained in $I$ and $JI^n = I^{n+1}$ for some nonnegative integer $n$.}
 of $I$.  Indeed, one readily
sees  that $JI = I^2$.   Moreover, the ideal $I$ is what Zariski calls a simple complete
ideal.\footnote{An ideal $I$ of an integral domain $R$ is said to be a {\it simple ideal} if $I$ is not the
unit ideal and I has no nontrivial factorization, that is, $ I =  JK$   implies either  $J$ or $K$
is the unit ideal of $R$. The ideal $I$ is {\it complete} if it is integrally closed. For integral
closure of ideals, see for example \cite{SH}.}

 The integral extension $ S[\frac{y}{x}]  = V$  can be seen to be  a DVR,  and  $V$ is the unique dicritical divisor
of $z = \frac{y^2}{x^3}$.  Let $v$ denote the valuation with value group $\mathbb Z$ associated
to the valuation ring $V$.  Then $v(y) = 3$ and $v(x) = 2$ and the image
$\overline z$ of $z$ in $V/M(V)$ is transcendental over $R/M$ and generates $V/M(V)$ over $R/M$.
\end{example}

  In  Definition~\ref{2.1}, we
define  the dicritical divisors  of a nonconstant bivariate polynomial.
%This definition relates directly  to the Jacobian problem.

\begin{definition}  \label{2.1}
Let $B$ denote the bivariate polynomial ring $k[X, Y]$ over a field $k$, and
let $L$ denote the quotient field of $B$. Let $I(B/k)$ denote  the set of DVRs $V$ on $L$
such that (i) $k \subset V$, (ii) the residue field of $V$ is transcendental over $k$, (iii)
$B \not\subset V$.

Let $f \in B \setminus k$ and let $I(B/k, f)$ denote the set of $V \in I(B/k)$ such that the
image of $f$ in the residue field of $V$ is transcendental over $k$.

Let
$B_f$ denote the localization of $B$ at the multiplicative set of nonzero elements of
$k[f]$. We observe that:
\begin{enumerate}
\item  $B_f = k(f)[X, Y]$   can be identified with 
the affine coordinate ring of  the generic curve $f^\sharp = 0$ where we take an indeterminate $u$ over $k$
and put $f^\sharp = f - u$.  This  identifies $B_f$ as the affine coordinate ring of $f^\sharp$ over 
the field $k(f)$;
\item the quotient field of $B_f$ is $L$; and
\item $\tr.deg_{k(f)} L = 1$.
\end{enumerate}
Consequently, $B_f$ is a one-dimensional
UFD and hence a PID. Therefore  the affine curve associated to $B_f$ is
irreducible and
nonsingular.  Let $I(B_f/k(f))$ be the set of all DVRs $V$  on $L$ that contain $k(f)$
and are such that $B_f \not\subset V$.
The elements in the  set $I(B_f/k(f))$
 are called the {\bf dicritical divisors  of} $f$  {\bf with respect to   the polynomial ring} $B$.
\end{definition}

\begin{remark} \label{1.2}
For $f \in B \setminus k$,   the set $I(B_f/k(f))$  of dicritical divisors of $f$ is  a finite set of DVRs.
It is equal to the set $I(B/k, f)$ of DVRs defined in Definition~\ref{2.1}.\footnote{Abhyankar elaborates
in \cite[pp.~147-156]{Ab9}  on this equality of the set $I(B/k, f)$ from surface theory with the set $I(B_f/k(f))$ from
curve theory.  The set $I(B/k, f)$ represents points at infinity of a  projective plane
curve,    while $I(B_f/k(f))$  is the set of all  
branches at infinity of  the  generic  curve   $f^\sharp$.}
 The dicritical divisors of $f$  are the places at infinity for the affine plane curve  having
coordinate ring $B_f$, cf. \cite[Remark~1, page~57]{AHS}.
\end{remark}

\begin{example}  \label{2.51}  Let $f = X^n   \in B$  with $n$ a positive integer.  Then the relative algebraic closure of
the field $k(f)$ in the field $L$ is the field $k(X)$.  It is straightforward to see that $f$ has one
dicritical divisor $V$   with respect to the polynomial ring $B$.  Moreover
$$
V ~:= ~ k(X)[Y^{-1}]_{Y^{-1}k(X)[Y^{-1}]}
$$
is that dicritical divisor.  If $n > 1$,  then $f$ is not a field generator.\footnote{An element $f \in B$ is a
{\it field generator} if there exists $\tau \in L$ such that   $f$ and $\tau$ generate $L$ as an
extension field of $k$, that is  $k(f, \tau) = L$. An element
 $f \in B$ is a {\it ring generator} if there exists $g \in B$  such that $k[f, g] = B$. }
\end{example}

\begin{example}  \label{2.515}  Let $f = X^mY^n   \in B$,  where  $m$ and $n$   are  positive integers
such that $\gcd(m, n) = 1$.
If  $a$ and $b$ are integers such that $mb - na = 1$   and $\tau$ is the rational function $X^aY^b \in L$,  then
$f$ and $\tau$ generate the field $L$ over $k$,  that is $k(f, \tau) = L$.  Hence $f$ is a field generator.
 It follows that  the field $k(f)$ is
relatively algebraically closed in the field $L$.   If $V$ is a dicritical divisor of the polynomial $f = X^mY^n$
with respect to the polynomial ring $B$, then $X^mY^n$ is a unit
of $V$  with the property that the image of $X^mY^n$ in  the residue field of $V$ is algebraically
independent over $k$.   Since $V$ does not contain $B$, either $X \not\in V$ or $Y \not\in V$.
 If $X \notin V$,  then $Y$ is in the maximal ideal of $V$ and vice versa.  We  conclude
 that  $X^mV = Y^{-n}V$, and the  polynomial $f$ has two dicritical divisors
with respect to the  polynomial ring $B$,  one that contains $X$ and another that contains $Y$.
It follows that  $f$ is not a ring generator, for if there exists an element
$g$ so that $k[f,g] = B$, the $f$ has only one dicritical divisor  $V$  with respect to
 ring $B$, namely $V = k(f)[g^{-1}]_{g^{-1}k(f)[g^{_1}]}$.
\end{example}

\begin{example}  \label{2.52}  Let $f = X^3 - Y^2$.  Then $f$ is irreducible in $B = k[X,Y]$ and the  affine
coordinate ring $\frac{B}{fB}$ of $f$ may be identified with the subring $k[t^2, t^3]$ of the polynomial ring $k[t]$
by means of the $k$-algebra homomorphism that maps $x \mapsto t^2$ and $y \mapsto t^3$.  Thus  $f$
defines a rational curve,  and the function field over $k$ of the affine coordinate ring $\frac{B}{fB}$ is a simple
transcendental field extension  $k(t)$.   Also the field $k(f)$ is relatively algebraically closed in $L$ and
$B_f = k(f)[X, Y]$   has one place at infinity.  Thus the polynomial $f$ has one dicritical divisor $V$ with
respect to the polynomial ring $B$.  The fractional ideal   $X^3V = Y^2V$ has  $V$-value $-6$,  and $f$ is a unit
in $V$ such that the image of $f$ is the residue field of $V$ is transcendental over the field $k$.
 The associated algebraic function field   $L/k(f)$ is
of  genus one.   Hence $L/k(f)$ is not a simple transcendental field extension.  The polynomial $f$
is not a field generator of $L/k$.
\end{example}

\begin{discussion}  \label{2.8}
We describe how the dicritical divisors of a bivariate polynomial are related to the dicritical divisors of an
element in the quotient field of
a two-dimensional regular local ring.   The polynomial ring $B = k[X, Y]$ is an affine component
 of  the   modelic projective plane   $\mathbb P^2_k   = \mathfrak M(k; X, Y, 1)$
 of the field $L$,   cf.  \cite[Section~5]{Ab10} or \cite[pages 116-119]{ZS2}.
Let  $\ell_{\infty}$ denote the line at infinity with respect to $B$.  Let
$$
f ~= ~ f(X,Y) ~= ~  \sum_{i + j \le N} a_{ij}X^iY^j ~ \in ~ k[X, Y]
$$
be a nonconstant polynomial, and let $V \in I(B_f/k(f))$.    Since $B \not\subset V$,  the
center of  $V$  on  $\mathbb P^2_k$ is a point on    $\ell_{\infty}$.  Let $R$ denote the two-dimensional regular local ring
associated to this point. Then $V \in \mathfrak D(R, f)$,  that is, $V$ is a dicritical divisor of $f$ in $R$.
Moreover,   $W \in  \mathfrak D(R, f)  \implies W \in I(B_f/k(f))$.   Thus the dicritical divisors in
 $ I(B_f/k(f))$ may be partitioned as follows: Let $R_1, \ldots, R_m$ be the  two-dimensional regular local
rings associated with points on $\ell_{\infty}$ that are the center of some $V \in I(B_f/k(f))$.  Then
 $ I(B_f/k(f))$  is the disjoint union of the sets $\mathfrak D(R_i, f)$.

Let $v$ denote a valuation associated to the valuation ring $V  \in   I(B_f/k(f))$.  Since $B \not\subset V$, either $v(X) < 0$
or $v(Y) < 0$.   Without loss of generality, we may assume $v(X) < 0$ and $v(X) \le v(Y)$.
Let $x  := 1/X$ and $y := Y/X$.  Then $R = k[x,y]_P$,  where $P = (x, \zeta(y))k[x, y]$ and
$\zeta(y) \in k[y]$ is an irreducible monic polynomial.  Let $\phi(X, Y) = \sum_{i + j = N} a_{ij}X^iY^j$  denote
the degree form of $f(X,Y)$.    Then $\phi(1,y) \in \zeta(y)k[y]$ and $x^Nf  \in R$.  Therefore $f$ generates a
special pencil.
\end{discussion}

\begin{definition}  \label{2.70}    Let $V \in \mathfrak D(R,z)$ be a dicritical divisor, and let $k'$ denote
the relative algebraic closure of the field  $k := R/M$ in  the residue field of $V$.
Then $k'/k$ is a finite algebraic field extension,   and
$V/M(V)$ is a simple transcendental extension $k'(\tau)$ of $k'$.   Moreover
the image $\overline z$ of
$z$ in $k'(\tau)$ is a nontrivial rational function in $\tau$.  The degree of the field extension $[k'(\tau): k(\overline z)]$ is
called the {\bf  degree of } $V$ {\bf  as a dicritical divisor in }  $\mathfrak D(R,z)$.
\end{definition}

In Example~\ref{1.3}, the degree of $V$ is one.  Example~\ref{2.7} describes a situation where
the degree of a dicritical divisor $V$ is greater than one.  Example~\ref{2.7}  also illustrates  how
 the local algebraic theory of dicritical divisors
 connects with that of polynomials in $\mathbb C[x,y]$.

\begin{example}   \label{2.7}  Consider the polynomial  $p(x,y) = x^4y^4 - x \in \mathbb C[x,y]$  and for each
$t \in \mathbb C$, let $C_t$ denote the
the curve
$ x^4y^4 - x - t  = 0$.   Notice that each point $(a,b) \in \mathbb C^2$ is on precisely one of the curves $C_t$ of the
pencil $\{C_t\}_{t \in \mathbb C}$.  Hence the polynomial $p$ defines a map $f_p:  \mathbb C^2 \to \mathbb C$.
Let  $\ell_{\infty}$ denote the line at infinity in $\mathbb P^2$  for
the natural embedding of $\mathbb C^2  ~  \hookrightarrow  ~ \mathbb P^2$ as in
Discussion~\ref{2.8}.  The curves $C_t$  have  two points on  the line $ \ell_{\infty}$,   the points $[1:0:0]$ and  $[0:1:0]$.
 Let $\widehat p(x,y,z) =  x^4y^4  -xz^7 -tz^8 \in \mathbb C[x,y,z]$  denote the homogenization of the polynomial $p - t$.
 We consider each of the two points on $\ell_{\infty}$:

 \noindent
 {\bf (1)}  For the point $[1:0:0]$,   set  $x = 1$ to get the rational function
 $t =  \frac{y^4 - z^7}{z^8}$.     Let   $R_1$  be the  $2$-dimensional regular local ring
 with maximal ideal  $ (y, z)R_1$ associated to this point.  We are interested in
 the dicritical divisors of the ideal $J_1 = ( y^4 - z^7, z^8)R_1$.  The integral closure of the ideal $J_1$ is  a simple
 complete ideal associated to a prime divisor $V$ with associated valuation $v$, where $v(y) = 7, v(z) = 4$ and $v(y^4 - z^7) = 32 = v(z^8)$.   Thus the ideal $J_1$ has
precisely one dicritical divisor and this dicritical divisor has degree one.
The integral closure
of $J_1$ is the simple complete ideal
$$
\overline{J_1} ~ = ~ (y^4 - z^7,~ z^8, ~ y^3z^3, ~ y^2z^5, ~ y^5)R_1.
$$

 \noindent
{\bf (2)}   For the point $[0:1:0]$,  set $y = 1$ to get the rational function $t  =  \frac{x^4 -  xz^7 }{z^8}$.   Let $R_2$ be
 the $2$-dimensional regular local ring  with maximal   ideal $(x,z)R_2$ associated to this point.
 We are interested in the dicritical divisors of the ideal
 $$
 J_2 ~=  ~(x^4 - xz^7,~ z^8)R_2.
 $$
  Consider the extension $R_2 \hookrightarrow R_2[\frac{x}{z}] :=  S$, and define $x_1:  =  \frac{x}{z}$.
  We have
$$
J_2S ~=~
(x^4 ~-xz^7, ~ z^8)S ~=~ (z^4x_1^4~- z^8x_1, ~ z^8)S ~
=~ z^4(x_1^4~-z^4x_1,~z^4)S.
$$
Thus the ideal $(x_1^4 - z^4x_1, z^4)S= (x_1^4, z^4)S$ is the
transform\footnote{If   $A \hookrightarrow B$ is a birational extension of unique factorzation domains and
$I$ is an ideal of $A$ not contained in any proper principal ideal, the {\it transform} of $I$ in $B$ is
the ideal $a^{-1}IB$,
where $aB$ is the smallest principal ideal in $B$ that contains $IB$ \cite[Definition~1.4]{L}.}
  of $J_2$ in $S$.
Since the transform of $J_2$ in $R_2[\frac{z}{x}]$ is the ring $R_2[\frac{z}{x}]$,  the ideal $J_2$ has a
unique base point in the blowup of the maximal ideal of $R_2$ with this unique base point being
the $2$-dimensional regular local ring $R_3 := S_{(z, x_1)S}$.
The
ideal $(x_1^4 - z^4x_1, z^4)R_3 = (x_1^4, z^4)R_3$ is the
 transform of $J_2$ in $R_3$.
The integral closure of $(x_1^4, z^4)R_3$  is $(x_1, z)^4R_3$,  the
4-th power of the maximal ideal of $R_3$. The integral closure of $J_2$
is the ideal
$$
\overline{J_2} ~=~ (x,~ z^2)^4R_2,
$$
and $(x,z^2)R_2$ is a simple complete ideal  associated to a prime
divisor $W$, with associated valuation $w$, where $w(x) = 2$ and
$w(z) = 1$.
 Thus the ideal $J_2$ has
precisely one dicritical divisor and this dicritical divisor has
degree four.
\end{example}

\begin{remark}   Let $R$ be a $2$-dimensional regular local ring and  let $J = (a,b)R$ be an ideal primary
for the maximal ideal of $R$.  Let $V$ be a dicritical divisor of $J$.  Let $k$ denote the residue field of $R$ and assume
that $k(\tau)$ is the residue field of $V$.  For $z = a/b$, let $\overline z$ denote the image of $z$ in
$V/M(V) = k(\tau)$.    It is natural to ask for conditions in order that   $\overline z$   be
a polynomial in $\tau$ as
opposed to just  being a rational function in $\tau$.
 Abhyankar and Luengo  in \cite{AbL}  prove this   for the prime
divisors at infinity of a nontrivial polynomial in a polynomial ring in 2 variables over a field.
More generally,  they prove in \cite{AbL} for $R$ a $2$-dimensional regular local domain that if $z$
generates a special pencil at $R$,  then $z$ generates a polynomial pencil in $R$.
\end{remark}

\section{Modelic spectrum and modelic blowup}

To describe a more general algebraic definition of dicritical divisors,  we review the following concepts.
Let $R$ be an integral domain.  The {\bf modelic spectrum } of $R$ is
$$
\frak V(R) ~ =  ~ \{ R_P ~| ~ P ~ \text{is a prime ideal of } ~ R \}.
$$
Here we  are identifying the prime ideals of $R$ with the local rings obtained by localizing
the integral domain  $R$ at these prime ideals.

\begin{definition}  \label{3.1}
Assume that $R$ is a subring of a field $K$ and let $x_0, \ldots, x_n$ be nonzero  elements of $K$.
The {\bf modelic proj} of $R$ with respect
to $x_0, \ldots, x_n$ is
$$
\frak W(R; x_0, \ldots, x_n) ~=~ \bigcup_{i=0}^n \frak V(R[ \frac{x_0}{x_i}, \ldots, \frac{ x_n}{x_i} ])
$$

Let $I$ be a nonzero ideal in the  integral domain $R$   and assume that  $x_0, \ldots, x_n$
are nonzero elements of $I$ that generate $I$.
 The {\bf modelic blowup } of $R$ at $I$ is
$$
\frak W(R, I)~= ~ \frak W(R; x_0, \ldots, x_n) ~=~ \bigcup_{i=0}^n \frak V(R[ \frac{x_0}{x_i}, \ldots, \frac{ x_n}{x_i} ])~=~ \bigcup_{0 \ne a \in I} ~\frak V(R[ I/a]).
$$
The modelic blowup of $R$ at $I$ is independent of the nonzero  generators chosen for  $I$.
Notice that $IR[I/a] = aR[I/a]$ is a principal ideal.

Let $t$ be an indeterminate over $R$ and let $R[It]$ denote the graded subring of the polynomial ring 
$R[t]$  generated by $It := \{at ~|~ a \in I \}$.   In other notation,   the modelic blowup 
of $R$ at $I$ is 
$$
\Proj R[It] ~: =~  \bigcup_{0 \ne a \in I} \Spec R[I/a].
$$
\end{definition}

\begin{discussion}  \label{3.2}
Consider the family of all quasilocal domains on the quotient field $K$ of $R$ that contain $R$, and
define a  partial order on this
set with respect to domination, i.e., for $S_1$ and $S_2$ define
$S_1 \le S_2$ if $S_1 \subset S_2$ and $M(S_2) \cap S_1 = M(S_1)$.

The quasilocal rings in the modelic blowup  $\frak W(R,I)$   of $R$ at $I$ are the minimal elements $S$
in  this partial order  such that  $IS$ is principal.
Every valuation domain  $V$  containing $R$ has a unique center on the modelic blowup
of $R$ at $I$,  i.e.,  there exists a unique quaslocal ring $S \in \frak W(R,I)$ such that $V$
dominates $S$.
\end{discussion}

Let $S$ be a quasilocal domain with quotient field $K$ and let $\overline S$ denote the
integral closure of $S$ in $K$.   Let
$S^{\mathfrak N}$ denote  the set of all members of the modelic spectrum $\mathfrak V(\overline S)$
of $\overline S$ that
dominate $S$.   The Lying Above Theorem for integral extensions \cite[T40, page~244]{Ab8}  implies
that the set  $S^{\mathfrak N}$ is nonempty,   and the Proper Containment Lemma  \cite[T42, page~245] {Ab8}
implies that
every quasilocal ring $T$ in the set
$S^{\mathfrak N}$ satisfies $\dim T \le \dim S$.
 If $\dim S$ is finite, then there
exists a quasilocal ring $T$ in $S^{\mathfrak N}$ with $\dim T = \dim S$  by the
Going Up Theorem \cite[T41, page~245]{Ab8}.
If $S$ is Noetherian, then $\dim S$ is finite by the Generalized Principal Ideal Theorem
\cite[T27, page~232]{Ab8}.   Moreover, by results proved by Mori and Nagata \cite[Theorem~33.10]{Nag},
the integral closure $\overline S$ of $S$ is a Krull domain and
the set $S^{\mathfrak N}$ is finite;  however, the quasilocal rings $T$ in $S^{\mathfrak N}$
may fail to be Noetherian.  This is because the integral closure
%$\overline{S}$
of a Noetherian local domain $S$ with
$\dim S \ge 3$ may fail to be Noetherian  \cite[Example~5, p. 207]{Nag}.

\begin{definition}  \label{3,3}
Let $I$ be a nonzero ideal in the  integral domain $R$   and assume that  $x_0, \ldots, x_n$
are nonzero elements of $I$ that generate $I$.
The {\bf normalized  modelic blowup } of $R$ at $I$ is
$$
\frak W(R, I)^{\mathfrak N} ~= ~ \frak W(R; x_0, \ldots, x_n)^{\mathfrak N} ~=~
 \bigcup_{i=0}^n \frak V(R[ \frac{x_0}{x_i}, \ldots, \frac{ x_n}{x_i} ])^{\mathfrak N}~=~
 \bigcup_{0 \ne a \in I} ~\frak V(R[ I/a])^{\mathfrak N}.
$$

\end{definition}

Concerning the modelic blowup  and normalized modelic blowup of ideals, it is natural to ask:

\begin{question} \label{3.4}  Let $I$ and $J$ be nonzero  finitely generated ideals of an integral domain $R$.
\begin{enumerate}
\item
Under what conditions does one have $\frak W(R,I) = \frak W(R,J)$?
\item
Under what conditions does one have $\frak W(R,I)^{\mathfrak N} = \frak W(R,J)^{\mathfrak N}$?
\end{enumerate}
Thus we are asking for conditions in order that the ideals $I$ and $J$ have the same blowup or the same
normalized blowup.
\end{question}

\begin{remark}  \label{3.5}  It is straightforward to see that an ideal $I$ and a power $I^n$ of $I$ have the
same modelic blowup and the same normalized modelic blowup.  Moreover,   if there exist nonzero elements $a$
and $b$ in $R$ such that $aI = bJ$,   then the ideals $I$ and $J$ have   the same modelic blowup and the same normalized modelic blowup.
\end{remark}

\begin{example} \label{3.6}
 Let $x$ and $y$ be indeterminates over a field $k$ and let $R = k[x^2, xy, y^2]$.
 The ring $R$ is the  coordinate ring of an affine surface that has an ordinary double point singularity at the
 origin.  Let $M  := (x^2, xy, y^2)R$ and let $P :=  (x^2, xy)R$.  Notice  that $M$ is a maximal ideal and $P$ is a
 height-one prime ideal.  Moreover,  we have $x^2M = P^2$.   Hence by Remark~\ref{3.5},  we have
 $\frak W(R, M) = \frak W(R, P)$.  The ideals $M$ and $P$ are both normal ideals.\footnote{An ideal $I$ of an integral domain
 is said to be a {\it normal ideal} if all the powers of $I$ are integrally closed.}   Therefore
 the modelic blowups of $M$ and $P$ are
 also their normalized modelic blowups and we have $\frak W(R, M)^{\mathfrak N} = \frak W(R, P)^{\mathfrak N}$.

 We also have $\frak W(R,M) =  \frak V(R[\frac{M}{x^2}])  \cup \frak V( R[\frac{M}{y^2}])$,
 and $R[\frac{M}{x^2}] = k[x^2, \frac{y}{x}]$  is isomorphic to a polynomial ring in 2 variables over the field $k$.
 Similarly,  we see that $R[\frac{M}{y^2}]$  is equal to $k[y^2, \frac{x}{y}]$.
 Therefore the modelic blowup of $R$ at $M$ is nonsingular.

\end{example}

\section{The Rees valuations and the dicritical divisors of an ideal}

Let $J$ be a nonzero ideal of a quasilocal domain $S$,  and  let
$\mathfrak W(S,J)^\Delta_1$ denote the set of all $1$-dimensional members $T$ of the blowup
$\mathfrak W(S,J)$ of $S$ at $J$ such that  $T$ dominates  $S$.

\begin{definition}  For a nonzero ideal $J$ of a quasilocal domain $S$, the set
$$
\mathfrak D(S,J) ~:=   ~(\mathfrak W(S,J)_1^\Delta)^{\mathfrak N}
$$ is the
{\bf dicritical set} of  $J$ in  $S$; members of
this set are called {\bf dicritical divisors} of $J$ in $S$.
\end{definition}

The set $\mathfrak D(S,J)$   may be empty.

For a nonzero ideal $I$ in a Noetherian integral domain $R$ one can take
the union of the sets  $\mathfrak{D}(R_P, IR_P)$,
where $P$ varies over the set of all prime ideals $P$ of $R$ that contain $I$.
This set may be described  as follows:

\begin{definition}
Let $I$ be a nonzero ideal in a Noetherian integral domain $R$, and let
$\mathfrak W(R,I,I)_1^\Delta$ denote all the one-dimensional members $S$ of the
blowup of $I$ such that $IS \ne S$.  Then $(\mathfrak W(R,I,I)_1^\Delta)^{\mathfrak N}$
is the set of {\bf dicritical divisors} of $I$.

\end{definition}

This gives a finite set of DVRs.

The blowup $\mathfrak W(R,I)$ is a  model that has finitely many one-dimensional
local rings $S$ that contain $I$.  For each such $S$, the set $S^{\mathfrak N}$ is
a finite set of DVRs.  The union over all $S$ gives
the dicritical divisors of $I$.

How does this compare with the set $\Rees I$ of Rees valuations rings of $I$?
The Rees valuation rings of $I$ are the DVRs that
contain $I$ and arise
as one-dimensional members of the normalized blowup $(\mathfrak W(R,I))^{\mathfrak N}$.

Thus $\Rees I = ((\mathfrak W(R,I,I))^{\mathfrak N})_1^\Delta$. The difference is one
first takes integral closure and then localizes.  For a nonzero ideal  $I$ of a Noetherian integral domain $R$,
the dicritical divisors of $I$ are Rees valuation rings of $I$, that is
$(\mathfrak W(R,I,I)_1^\Delta)^{\mathfrak N}$  is  a subset of $\Rees I$.
 If $R$ is also
 universally catenary,  then  the
Rees valuations of $I$ are precisely the same as the dicritical divisors of $I$.
Thus if $I$ is an $\m$-primary ideal of a universally catenary Noetherian local domain $(R,\m)$ then
the set $\Rees I$ of Rees valuation rings of $I$ is precisely the set  $\mathfrak D(R,I)$  of dicritical divisors
of $I$.

The following three results  about Rees valuation rings  of an ideal  are given in \cite{HKT}.

\begin{theorem} \label{3.1hkt}
Let $(R,\m)$ be a universally catenary  analytically
unramified Noetherian local domain with $\dim R = d$, and
let $V$ be a prime divisor of $R$ centered on $\m$. Let $I \subseteq \m$ be
an ideal of $R$. The following are equivalent
\begin{enumerate}
\item $V \in \Rees I$.
\item There exist elements $b_1, \ldots, b_d$ in $I$ such that
$b_1V = \cdots= b_dV = IV$
and the images of  $~\frac{b_2}{b_1}, \ldots, \frac{b_d}{b_1}$ in the residue
field $k_v$ of $V$ are algebraically
independent over $R/\m$.
\item If $I = (a_1, \ldots, a_n)R$, then there exist
elements $b_1, \ldots, b_d$ in $\{a_i\}_{i=1}^n$ such that
$b_1V = \cdots= b_dV = IV$
and the images of $~\frac{b_2}{b_1}, \ldots, \frac{b_d}{b_1}$ in the residue
field $k_v$ of $V$ are algebraically
independent over $R/\m$.
\end{enumerate}
Thus  if $I = (a_1, a_2, \ldots, a_d)R$, then $V \in \Rees I  \iff  a_1V = a_2V =\cdots= a_dV$
and the images of $\frac{a_2}{a_1}, \ldots, \frac{a_d}{a_1}$ in
$k_v$ are algebraically
independent over $R/\m$.
\end{theorem}

\begin{proposition}\label{contracted}
Let $(R, \m)$ be a $d$-dimensional regular local ring with $d \ge 2$, let $x \in \m \setminus \m^2$, and let $S = R [\frac{\m}{x}]$.
Let $K$ be an ideal of $R$ that is contracted from $S$, and let $x f \in K$, where $f \in R$.
Then
\begin{enumerate}
\item
$g f \in K$ for each $g \in \m$.
\item
If $\ord_{R} (x f) = \ord_{R} (K)$, then the order valuation $\ord_R$
is a Rees valuation of $K$.
\end{enumerate}
\end{proposition}

\begin{proposition} \label{reestrans}
Let $(R,\m)$ be a $d$-dimensional regular local ring and
let $(S,\n)$ be a $d$-dimensional regular local ring that birationally dominates $R$.
Let $I$ be an $\m$-primary ideal of $R$ such that its transform $J = I^{S}$  in $S$
is not equal to $S$, and let $V$ be a DVR that  birationally dominates $S$.
Then
$$
V \in ~ \Rees_S J   ~ \iff   ~ V \in ~ \Rees_R I.
$$
\end{proposition}

\begin{remark} With notation as in Theorem~\ref{3.1hkt}, let $\overline{\frac{b_2}{b_1}}, \ldots,
\overline{\frac{b_d}{b_1}}$ denote the images of $\frac{b_2}{b_1}, \ldots, \frac{b_d}{b_1}$ in
the residue field $k_v$ of $V$. An interesting integer associated with $V \in \Rees I$ and $b_1, \ldots b_d$ is the field degree
$$
\Big [k_v: (R/\m)\Big (\overline{\big (\frac{b_2}{b_1}\big )},\ldots,  \overline{\big (\frac{b_d}{b_1}\big )}\Big )\Big ]
$$
Notice the analogy with the degree of $V$ as a dicritical divisor in Definition~\ref{2.70}.
\end{remark}

\section{Dicritical divisors and the structure of quadratic sequences}

A central question concerning dicritical divisors that interested Abhyankar can be stated as follows:

\begin{question}  \label{5.1}   Let $R$ be a 2-dimensional regular local domain with quotient field $L$
and let $U \subset D(R)^{\Delta}$ be a finite set of prime divisors of $R$.
\begin{enumerate}
\item Does there exist an
element $z \in L$ such that $\mathfrak D(R,z) = U$?
\item  If the answer to item~1  is affirmative, is it possible in
some algorithmic way to describe all the $z \in L$ such that $\mathfrak D(R,z) = U$?
\end{enumerate}
\end{question}

In a sequence of papers written with Heinzer \cite{AH1}, \cite{AH2}, \cite{AH3},  Question~\ref{5.1} and
various related questions are studied. In particular,  an affirmative answer to item~1 of Question~\ref{5.1} is
given in \cite{AH1}.
If the residue field of $R$ is infinite,  it is also shown in \cite{AH1} that the element $z$ can be chosen
so that,  in the terminology of Definition~\ref{2.70},  for each $V \in U$,   the   degree of  $V$   
as a dicritical divisor in   $\mathfrak D(R,z)$  is 1.

A main tool  in \cite{AH1} is the work  of Zariski   in Appendix 5 of \cite{ZS2} concerning the structure of complete
ideals of a $2$-dimensional regular local domain.  Let $Q(R)$ denote the set of 2-dimensional regular local
domains that birationally dominate $R$.  For each $S \in Q(R)$, the order valuation domain $\ord S$ is a
prime divisor on $R$, and  Theorem~\ref{1.1}  implies that the map $Q(R) \to D(R)^{\Delta}$ that maps
$S$ to $\ord S$ is a bijection of the sets $Q(R)$ and $D(R)^{\Delta}$.

Let $I$ be a simple complete $M$-primary ideal.  If $I \ne M$,  Zariski proves the existence of
  a positive integer $\nu$ and a  unique regular local ring $S \in Q(R)$  such that the
finite sequence  as given in Theorem~\ref{1.1}:
$$
R~= ~ R_0 ~\subset ~R_1 ~\subset ~ \cdots ~ \subset R_{\nu} ~= ~ S,
$$
where $R_i$ is a local quadratic transform  of $R_{i-1}$ for each $i \in \{1, \ldots, \nu\}$,
consists precisely of the regular local rings $T \in Q(R) $ for which the transform of $I$ in $T$ is
a proper ideal of $T$.   Moreover,  the transform of $I$ in $S$ is the maximal ideal of $S$.
The regular local domains  $R_i$,   with   $i \in \{1, \ldots, \nu\}$,   are the base points
of $I$.\footnote{If $I$ is an $M$-primary ideal of a regular local
ring $R$,  the {\it base points} of $I$ are the regular local rings $S$ infinitely near to $R$ such that
the transform of $I$ in $S$ is not the unit ideal.  The ideal $I$ is said to be {\it finitely supported} if it 
has only finitely many base points.  If $\dim R = 2$,  then every $M$-primary ideal is finitely supported. This 
is no longer true if $\dim R \ge 3$.}

Let
$C(R)$ denote the set of $M$-primary simple
complete ideals of  $R$.   The Zariski Quadratic Theorem \cite[page~173]{AH1} asserts that
for each $ V \in D(R)^{\Delta}$ there exists  at least one and at most a finite number of
$V$-ideals\footnote{An ideal $J$ of $R$ is a $V$-{\it ideal} if $J = JV \cap R$.}
in $R$ that are members of $C(R)$.  Labeling these $V$-ideals of $C(R)$ as
$$
M~ = ~J_0 ~ \supset  ~ J_1 ~ \supset ~ \cdots ~\supset J_{\nu}
$$
one obtains a bijection $\zeta_R:  D(R)^{\Delta} \to C(R)$ by defining $\zeta_R(V) = J_{\nu}$.

Zariski proves that a product of complete ideals in $R$ is again complete and
every complete $M$-primary ideal  can be expressed uniquely as
a finite product of powers of the ideals in $C(R)$.    For a finite subset $U$ of $D(R)^{\Delta}$,
define $\overline{\zeta}_R(U)$ to be the product $\prod_{V \in U}\zeta_R(V)$.
Using results from Northcott-Rees \cite{NoR},  it is shown in \cite{AH1} that some power of the ideal
$\overline{\zeta}_R(U)$ has a 2-generated reduction.  This yields an affirmative answer to item~1 of
Question~\ref{5.1}.

To consider item~2 of Question~\ref{5.1},  Abhyankar turned to a careful and detailed analysis of
quadratic transformations (QDT's) and inverse QDT's.  Inverse transforms as considered by
Zariski \cite[pp. 390-391]{ZS2} involve  complete ideals.  Abhyankar's interest was more in
2-generated reductions.  For a complete  ideal   $I$  that  is  primary to the maximal ideal
 of a $2$-dimensional regular local domain,  the ideals $J$ that are 2-generated reductions of
 $I$ are usually very far from being unique.\footnote{The intersection of all the reductions of $I$ is
 called the {\it core}  of $I$.  Interesting results about he core of  a complete ideal $I$
 that is primary for the maximal ideal of a $2$-dimensional regular local domain with infinite
 residue field are obtain by Huneke and
 Swanson in \cite{HS}. They prove that the core of $I$ is integrally closed and is the the product of $I$
 with its second Fitting ideal.}
Abhyankar  was more interested in specific properties of QDT's and
their inverses.
 In this connection  he coined the term antiquadratic
transformation sequence.  Example~\ref{5.2} below, illustrates the type  of
 results  Abhyankar proved  to  answer  item~2 of Question~\ref{5.1} in special cases.

To construct  examples such as Example~\ref{5.2},  and to develop the algebraic theory of curvettes  as
in \cite{AbA},  Abhyankar developed another characterization of dicritical divisors.  This characterization
depends on the concept  of Zariski number and Zariski index  that he defined as follows:

\begin{definition} \label{5.15}
Let $R$ be a $2$-dimensional regular local ring with maximal ideal $M$ and let 
$\grad R  = R/M \oplus M/M^2 \oplus \cdots $ denote the 
associated graded ring of $R$.  The ring $\grad R$ is a polynomial ring in 2 variables over the 
residue field $k := R/M$.   If $M = (x,y)R$, then the initial forms $\info x$ and $\info y$ of $x$ and $y$ in 
$M/M^2$ are algebraically independent over $k$ and generate $\grad R$.  Let $J$ be a nonzero 
ideal in $R$ and let $d := \ord_R J$. Thus $J \subset M^d$ and $J \not\subset M^{d+1}$. 
Let $s$ denote the degree of the GCD of all the initial forms  $\info f$ of 
 elements  $f \in J$  such that $\ord_R f = d$. The {\bf Zariski number} $m(R,J)$ is defined by $m(R,J):=d-s$. 
 Notice that the Zariski number of $J$ is a nonnegative integer less than or equal to $d$.  If $J$ is 
 $M$-primary, Zariski proved that the power of $M$ that occurs in the factorization of the integral 
 closure of $J$ as a product of simple complete ideals is the integer $m(R,J)$.  
 Let $T\in Q(R)$ be a quadratic transform of $R$.   The 
{\bf derived Zariski number} $m(R,J,T)$ is  the nonnegative integer 
 $m(R,J,T):=m(T,J_T)$,  where $J_T$ denotes the transform of  the ideal $J$ in $T$. 
 For a prime divisor  $V\in D(R)^\Delta$, 
the {\bf Zariski index} $n(R,J,V)$ is  the nonnegative integer
  $n(R,J,V)=m(R,J,o_R^{-1}(V))$, where $o_R^{-1}(V)$ is the unique $2$-dimensional regular local 
  domain $S \in Q(R)$ such that $\ord_S = V$.    Using this notation,  the 
  Zariski Factorization Theorem concerning complete ideals of the $2$-dimensional regular
  local domain $R$ may be stated as follows:
  \begin{equation*}
  J^{-R} ~ ~ =    ~\text{GCD}(J)_R\prod_{V\in D(R)^\Delta}\zeta_R(V)^{n(R,J,V)},  
  \end{equation*}
where
  $J^{-R}$ denotes the integral closure of $J$  and $\text{GCD}(J)_R$ denotes the smallest nonzero principal ideal in $R$ containing $J$.  The set   $\mathfrak D(R,J)$ of dicritical divisors of $J$ is 
  characterized as 
  $$
  \mathfrak D(R,J)  ~=   ~  \{V\in D(R)^\Delta  ~|~ n(R,J,V)>0\}.
  $$
    A prime divisor $V$ is a dicritical
  divisor of the ideal $J$ if and only if the Zariski index $n(R,J,V)$ is positive.
\end{definition}  

In  
Example~\ref{5.16}  we illustrate the Zariski number,  derived Zariski number and Zariski index in a 
simple example.

\begin{example}  \label{5.16}
Let $R$ be a  $2$-dimensional regular local domain with maximal ideal $M(R) = M = (x,y)R$,
and let $J := (x^3, x^2y, y^7)R$.  We have $\ord_R J = 3$,  and the GCD of $\info x^3$ and $\info x^2y$ 
is a polynomial in $\grad R$ of degree 2. Hence the Zariski number $m(R, J) = 3 - 2 = 1.$   Thus $M$ 
divides the integral closure of $J$, but $M^2$ does not.  The ideal $J$ has one base point in the first
neighborhood of the blowup of $M$, namely  $R_1 := R[\frac{x}{y}]_{(y, \frac{x}{y})R[\frac{x}{y}]}$.  
Let $x_1 := \frac{x}{y}$.  We have $JR_1 = (y^3x_1^3, y^3x_1^2, y^7)R_1$.   The transform of $J$ in 
$R_1$ is the ideal $J_1 := (x_1^3, x_1^2, y^4)R_1$.   We have $\ord_{R_1}J_1 = 2$ and the 
derived Zariski number $m(R, J, R_1) = 0$.   Thus the maximal ideal $M_1$ of $R_1$ does not 
divide the integral closure of $J_1$ in $R_1$.  Let 
$R_2 := R_1[\frac{x_1}{y}]_{(y, \frac{x_1}{y})R_1[\frac{x_1}{y}]}$.  It is straighforward to check that $R_2$ 
is the only base point of $J$ in the second neighborhood of $R$.  Let $x_2 := \frac{x_1}{y}$.  We have
$J_1R_2 = (y^3x_2^3, y^2x_2^2, y^4)R_2 = y^2(yx_2^3, x_2^2, y^2)R_2$.  The transform of $J$ in $R_2$ 
is the ideal $J_2  := (x_2^2, y^2)R_2$.  We have $\ord_{R_2} J_2 = 2$ and the derived Zariski number
$m(R, J, R_2) = 2$.   It follows that the integral closure of $J_2$ is $M_2^2$,  where $M_2$ is the 
maximal ideal of $R_2$.  The ideal $J$ has 3 base points $R = R_0, R_1, R_2$.  Let $V_i := \ord_{R_i}, 
i  \in  \{0,1,2 \}$.  The Zariski index $n(R,J, V_0) = 1$,  while $n(R,J, V_1) = 0$ and $n(R,J, V_2) = 2$.  
The ideal $J$ has two dicritical divisors, namely $V_0$ and $V_2$.   
\end{example} 

\begin{example}   \label{5.2}
Let $R$ be a  $2$-dimensional regular local domain with maximal ideal $M(R) = M = (x,y)R$.
Consider the infinite QDT-sequence $(S_j)_{0\le j<\infty}$   such that  $S_0 = R$ and
$M(S_j)=(x_j,y_j)S_j$, where $x = x_j$ and $y_j = \frac{y}{x^j}$ for each $j \ge 1$.  Let
$V_j := \ord S_j$.   The simple complete ideal $\zeta_R(V_j) = (x, y^j)R$ for each $j$.   For each
$m \in \N$,  let $I_m=\prod_{0\le j<m}\zeta_R(V_j)$.  The complete ideal $I_m$ has order $m$ and is
minimally generated by $m+1$ monomials in $x$ and $y$.  Let $B(n) = \frac{n(n+1)}{2}$   and consider
the polynomials
$$
 F_m(X,Y) :=\sum_{0\le p\le[\frac{m-1}{2}]}X^{B(2p+1)}Y^{m-1-2p} ~~\text{,}~ ~
 G_m(X,Y)=\sum_{0\le p\le[\frac{m}{2}]}X^{B(2p)}Y^{m-2p}
$$
over the ring of integers.  Abhyankar showed the ideal $J_m := (F_m(x,y), G_m(x,y))R$ is a $2$-generated
reduction of $I_m$ for each $m \in \N$.
Thus,  for example, if
$m=5$, then the  ideal
$$
J_5 ~=  ~ (xy^4+x^6y^2+x^{15},  ~y^5+x^3y^3+x^{10}y)R
$$
is a 2-generated reduction of the complete ideal $I_5$.  

 Notice that  for each prime divisor $V \in  D(R)^\Delta$ and each $m \in \N$,  
 the  Zariski index  $n(R,J_m,V)$   of the ideal $J_m$ is
 either $0$ or $1$.  Moreover,  $n(R,J_m,V)  = 1$ if and only if $o_R^{-1}(V) \in \{S_0,\dots,S_m\}$. 
 \end{example}

 An attractive feature of Example~\ref{5.2} 
 is
that it applies  without any restrictions on the residue field of $R$.
The residue field of $R$ could be the finite field with two elements.

\begin{remark}
Let $(R, M)$ be a Noetherian local ring with $\dim R = d > 0$ and let $I$ be an $M$-primary ideal.  If the
residue field of $R$ is infinite, then Northcott and Rees \cite{NoR}   prove  that there exist
 $d$-generated ideals $J$
that are reductions of $I$.  Moreover, each $d$-generated reduction $J$ of $I$ is a minimal reduction
in the sense that there is no ideal properly contained in $J$ that is a reduction of $I$.  If $R$ has a finite
residue field,  it may happen that $I$ fails to have a $d$-generated reduction.  Let $F$ be an
arbitrary finite field.
 In \cite[Example~2.3]{HRR},  an example is given of a $2$-dimensional
 Cohen-Macaulay local ring $(R,M)$ such that
 $R$ has residue field $F$ and the maximal ideal $M$ of $R$ fails to have a $2$-generated reduction.
\end{remark}

Specific solutions to  item~2 of Question~\ref{5.1}  such at those obtained by
Abhyankar in Example~\ref{5.2} and those obtained by
Abhyankar and Artal in \cite[Theorem~3.6 and Theorem~4.5]{AbA}  indicate
that Question~\ref{5.4} may possibly have an affirmative answer.

\begin{question}  \label{5.4}    Let $(R,M)$ be a 2-dimensional regular local domain having a finite residue field,
and let $I$ be a complete  $M$-primary ideal.  Does there  always  exist a $2$-generated ideal
$J = (a,b)R$ such that $J$ is a reduction of $I$?
\end{question}

\bibliographystyle{amsalpha}

\end{document}